\documentclass[10pt]{amsart}

\makeatletter
\usepackage{amssymb}
\usepackage{latexsym}
\usepackage{amsbsy}
\usepackage{amsfonts}

\def\marginpar#1{\ignorespaces}

\newtheorem{theorem}[equation]{Theorem}
\newtheorem{proposition}[equation]{Proposition}
\newtheorem{lemma}[equation]{Lemma}
\newtheorem{corollary}[equation]{Corollary}
\newtheorem{definition}[equation]{Definition}

\theoremstyle{definition}
\newtheorem{remark}[equation]{Remark}

\newtheorem{example}[equation]{Example}

\numberwithin{equation}{section}

\def\AArm{\fam0 \rm}%
\newdimen\AAdi%
\newbox\AAbo%
\def\AAk#1#2{\setbox\AAbo=\hbox{#2}\AAdi=\wd\AAbo\kern#1\AAdi{}}%

\newcommand{\BBone}{{\ensuremath{{\AArm 1\AAk{-.8}{I}I}}}}

\def\eqref#1{(\ref{#1})}
\def\eqlabel#1{\def\@currentlabel{#1}}

\def\formula#1{\def\@tempa{#1}\let\@tempb\theequation\def\theequation{%
\hbox{#1}}\def\@currentlabel{(\theequation)}$$}
\def\endformula{\leqno\hbox{(\@tempa)}$$\@ignoretrue\let\theequation\@tempb}

\def\given{\hskip5\p@\relax\vrule\@width.4\p@\hskip5\p@\relax}

\newcommand{\open}[1]{%
\par\normalfont\topsep6\p@\@plus6\p@\trivlist\item[\hskip\labelsep\itshape#1%
\@addpunct{.}]\ignorespaces}

\DeclareRobustCommand{\close}[1]{%
  \ifmmode 
  \else \leavevmode\unskip\penalty9999 \hbox{}\nobreak\hfill
  \fi
  \quad\hbox{$#1$}}

\newlength{\toskip}

\def \R {{\mathbb R}}

\def \P {{\mathbb P}}
\def \E {{\mathbb E}}

\def \D {{\mathbb D}}
\def \L {{\mathbb L}}

\def \phi {\varphi}

\makeatother

\begin{document}
\date{\today}

\title[HYPERCONTRACTIVITY ...]{HYPERCONTRACTIVITY FOR PERTURBED DIFFUSION SEMIGROUPS.}

 \author[P. Cattiaux]{\quad {Patrick} Cattiaux }

\address{{\bf {Patrick} CATTIAUX},\\ Ecole Polytechnique, CMAP, F- 91128 Palaiseau cedex,
CNRS 756\\ and Universit\'e Paris X Nanterre, \'equipe MODAL'X, UFR SEGMI\\ 200 avenue de
la R\'epublique, F- 92001 Nanterre, Cedex.} \email{cattiaux@cmapx.polytechnique.fr}

\maketitle

 \begin{center}
 \textsc{Ecole Polytechnique \quad and \quad Universit\'e Paris X}
 \end{center}

\begin{abstract}
$\mu$ being a nonnegative measure satisfying some Log-Sobolev inequality, we give conditions on $F$
for the Boltzmann measure $\nu \, = \, e^{-2F} \mu$ to also satisfy some Log-Sobolev inequality.
This paper improves and completes the final section in \cite{CX02}. A general sufficient condition
and a general necessary condition are given and examples are explicitly studied.\\
\par\vspace{10pt}
\noindent{\sc R\'esum\'e}. $\mu$ \'etant une mesure positive satisfaisant une in\'egalit\'e de
Sobolev logarithmique, nous donnons des conditions sur $F$ pour que la mesure de Boltzmann $\nu \, =
\, e^{-2F} \mu$ satisfasse \'egalement une telle in\'egalit\'e (am\'eliorant et compl\'etant ainsi
la derni\`ere partie de \cite{CX02}). Les conditions obtenues sont illustr\'ees par des exemples.
\end{abstract}
\bigskip

\bigskip

\bigskip

\section{\bf Introduction and Framework.}\label{I}

In \cite{CX02} we have introduced a pathwise point of view in the study of classical inequalities.
The last two sections of this paper were devoted to the transmission of Log-Sobolev and Spectral Gap
inequalities to perturbed measures, without any explicit example. In the present paper we shall
improve the results of section 8 in \cite{CX02} and study explicit examples. Except for one point,
the present paper is nevertheless self-contained. In order to describe the contents of the paper we
have first to describe the framework.
\medskip

{\bf Framework.} For a nonnegative measure $\mu$ on some measurable space E, let us first consider a
$\mu$ symmetric diffusion process $(\P_x)_{x\in E}$ and its associated semi-group $(P_t)_{t \geq 0}$
with generator $A$. Here by a diffusion process we mean a strong Markov family of probability
measures $(\P_x)_{x\in E}$ defined on the space of continuous paths $\mathcal C^0(\R^+,E)$ for some,
say Polish, state space $E$, such that there exists some \underline{algebra} $\D$ of uniformly
continuous and bounded functions (containing constant functions) which is a core for the extended
domain $D_e(A)$ of the generator (see \cite{CL96}).

One can then show that there exists a countable orthogonal family $(C^n)$ of local martingales and a
countable family $(\nabla^n)$ of operators s.t. for all $f\in D_e(A)$
\begin{equation}\label{1.1}
M_t^f=f(X_t)-f(X_0)-\int_0^t \, Af(X_s) \, ds=\sum_n \, \int_0^t \, \nabla^n f(X_s)
\, dC^n_s \, ,
\end{equation}
in $\mathbb M^2_{loc}(\P_{\eta})$ (local martingales) for all probability measures $\eta$ on $E$.

One can thus define the ``carr\'e du champ'' $\Gamma$ by $\Gamma(f,g)=\sum_n \, \nabla^n f \,
\nabla^n g  \stackrel{\rm def}{=} (\nabla f)^2$ , so that the martingale bracket is given by
$<M^f>_t \, =\int_0^t \, \Gamma(f,f)(X_s) \, ds \, .$ In terms of Dirichlet forms, all this, in the
symmetric case, is roughly equivalent to the fact that the local pre-Dirichlet form $\mathcal E
(f,g)=\int \, \Gamma(f,g) \, d\mu \quad f,g\in \D$ is closable, and has a regular (or quasi-regular)
closure $(\mathcal E, D(\mathcal E))$, to which the semigroup $P_t$ is associated. Notice that with
our definitions, for $f\in \D$
\begin{equation}\label{1.2}
\mathcal E(f,f)=\int \, \Gamma(f,f) \, d\mu=-2 \, \int \, f \, Af \, d\mu \, =
 \,- \, \frac{d}{dt} \, \parallel P_tf\parallel^2_{\L^2(\mu)}|_{t=0} \, .
\end{equation}
It is then easy to check that $\Gamma(f,g)=A \, (fg) - f \, Ag -g \, Af$ , that $\D$ is stable for
the composition with compactly supported smooth functions and satisfies the usual chain rule.
\bigskip

{\bf Content.} The aim of this paper is to give conditions on $F$ for the perturbed measure $\nu_F
\, = \, e^{-2F} \mu$ to satisfy some Logarithmic Sobolev inequality, assuming that $\mu$ does. As in
the final section of \cite{CX02} these conditions are first described in terms of some martingale
properties in the spirit of the work by Kavian, Kerkyacharian and Roynette (see \cite{KKR}) (see
section \ref{III} Theorem \ref{3.10}).
\smallskip

We shall then study in section \ref{IV} how this general criterion can be checked in the same
general situation . Here again we are inspired by \cite{KKR} (Well Method). It turns out that the
Well Method can be generalized to other $F$-Sobolev inequalities (see \cite{BCR}).
\medskip

Since sections \ref{III} and \ref{IV} are concerned with the hyperbounded point of view, and
following the suggestion of an anonymous referee, we study in section \ref{ls} the log-Sobolev point
of view (i.e. the perturbation point of view is analyzed on log-Sobolev inequalities). We show that
both point of view yield (almost) the same results.
\medskip

In the final section we study some examples, namely Boltzmann measures on $\R^N$. Explicit examples
and counter examples are given, and some comparison with existing results is done.
\bigskip

\textbf{Acknowledgements.} I wish to thank Michel Ledoux for his interest in this work and for
pointing out to me Wang's results. I also benefited of nice discussions with Franck Barthe, Cyril
Roberto and Li Ming Wu.
\bigskip

\begin{center}
{\bf  Some notation and general results}
\end{center}

The material below can be found in many very good textbooks or courses see e.g. \cite{Ane},
\cite{Bak94}, \cite{DSbook}, \cite{Gross}, \cite{GZ99}, \cite{Led99}.
\medskip

We shall say that $\mu$ satisfies a Log-Sobolev inequality LSI if for some universal constants $a$
and $b$ and all $f\in \D \, \cap \, L^1(\mu)$,
\begin{equation}\label{2.1}
\int \, f^2 \, \log \Big(\frac{f^2}{\parallel f\parallel_{\L^2(\mu)}^2}\Big) \, d\mu \, \leq \, a \,
\int \, \Gamma(f,f) \, d\mu \, + b \, \parallel f\parallel_{\L^2(\mu)}^2 \, .
\end{equation}

When $b=0$ we will say that the inequality is tight (TLSI), when $b>0$ we will say that the
inequality is defective (DLSI). So we never will use (LSI) without specifying (TLSI) or (DLSI).

Note that when $\mu$ is bounded \eqref{2.1} easily extends to any $f \in D(\mathcal E)$. It is not
the case when $\mu$ is not bounded, in which case it only extends to $f \in D(\mathcal E) \, \cap \,
\L^1(\mu)$ or to $f \in D(\mathcal E)$ but replacing $\log$ by $\log^+$ in the left hand side of
\eqref{2.1}. An example of such phenomenon is $f=(1+|x|)^{- \, \frac 12} \, \log^{\alpha}(e+|x|)$
for $1< \, 2 \alpha < \, 2$, $E=\mathbb R$ and $d\mu=dx$.
\smallskip

 These inequalities are known to be related to continuity or contractivity of the semigroup
$P_t$. We shall say that the semigroup is hyperbounded (resp. hypercontractive) if for some $t>0$
and $p>2$, $P_t$ maps continuously $\L^2(\mu)$ into $\L^p(\mu)$ (resp. is a contraction). In this
case we shall denote the corresponding norm $\parallel P_t\parallel_{\L^2(\mu) \, \rightarrow \,
\L^p(\mu)} \, $, or simply $\parallel P_t\parallel_{2,p}$ when no confusion is possible. It is well
known that hyperboundedness (resp. hypercontractivity) is equivalent to (DLSI) (resp. (TLSI)) (see
e.g. \cite{Bak94} Theorem 3.6 or \cite{CX02} Corollary 2.8). Gross theorem tells next that
boundedness or contraction hold for all $p>2$ for some large enough $t$. Replacing $p$ by $+\infty$
in the definition we get the notion of ultracontractivity extensively studied in the book by E.B.
Davies \cite{Dav}. Links with Log-Sobolev inequalities are especially studied in chapter 2 of
\cite{Dav}.

Finally recall that (TLSI) is equivalent to (DLSI) plus some spectral gap condition (as soon as we
will use spectral gap properties we shall assume that $\mu$ is a probability measure). The usual
spectral gap (or Poincar\'e) inequality will be denoted by (SGP). A weaker one introduced by
R{\"o}ckner and Wang (see \cite{RW01}) called the weak spectral gap property (WSGP) is discussed in
\cite{Aid01} and in section 5 of \cite{CX02}. In particular (DLSI)+(WSGP) implies (TLSI) originally
due to Mathieu (\cite{Mat98}) is  shown in \cite{CX02} Proposition 5.13.
\bigskip

\section{\bf Hypercontractivity for general Boltzmann measures.}\label{III}

We introduce in this section a general perturbation theory. In the framework of section \ref{I} let
$F$ be some real valued function defined on $E$.
\begin{definition}\label{3.1}
The Boltzmann measure associated with $F$ is defined as $\nu_F \, = \, e^{-2F} \, \mu$.
\end{definition}
When no confusion is possible we may not write the subscript $F$ and simply write $\nu$.
\medskip

The transmission of Log-Sobolev or Spectral Gap inequalities to Boltzmann measures has been
extensively studied in various contexts. The first classical result goes back to Holley and Stroock.

\begin{proposition}\label{3.2}
Assume that $\mu$ is a probability measure and $F$ is bounded. Then if $\mu$ satisfies (DLSI) with
constants $(a,b)$, $\nu_F$ satisfies (DLSI) with constants $(a \, e^{Osc(F)} \, , \, b \,
e^{Osc(F)})$ where $Osc(F)=\sup(F) \, - \, \inf(F)$.
\end{proposition}
This result is often stated with $2 \, Osc(F)$ i.e. with a useless factor 2 (see \cite{Ro99}
Proposition 3.1.18).
\smallskip

When $F$ is no more bounded, general (though too restrictive) results have been shown by Aida and
Shigekawa \cite{AS94} (also see \cite{CX02} section 7). Other results can be obtained through the
celebrated Bakry-Emery criterion. As in section 8 of \cite{CX02} we shall follow a beautiful idea of
Kavian, Kerkyacharian and Roynette (see \cite{KKR}) in order to get better results (with a little
bit more regularity). The main idea in \cite{KKR} is that ultracontractivity for a Boltzmann measure
built on $\R^N$ with $\mu$ the Lebesgue measure and $F$ regular enough, reduces to check the
boundedness of \underline{one} and \underline{only one} function.

The aim of this section is to improve these results. First let us state the hypotheses we need for
$F$.

\begin{enumerate}\stepcounter{equation}\eqlabel{\theequation}\label{3.3}
\item[\theequation] \textbf{Assumptions H(F)}
\smallskip

\item[(1)] \quad $\nu_F$ is a probability measure, $F\in D(\mathcal E)$ ,
\item[(2)] \quad for all $f\in \D$, $\mathcal E_F(f,f)=\int \, \Gamma(f,f) \, d\nu_F < +\infty$ ,
\item[(3)] \quad for all $f\in \D$, $Af \in \L^1(\nu_F)$ ,
\item[(4)] \quad $\int \, \Gamma(F,F) \,  d\nu_F < +\infty$ .
\end{enumerate}

The Girsanov martingale $Z_t^F$ is then defined as
\begin{eqnarray}\label{3.4}
Z_t^F & = & \exp \, \{- \, \int_0^t \, \nabla F(X_s).dC_s - \, \frac 12 \, \int_0^t \, \Gamma(F,F)
(X_s) \,  ds\} \, .
\end{eqnarray}
When H(F) holds, we know that $Z_.^F$  is a $\P_x$ martingale for $\nu_F$, hence $\mu$ almost all
$x$. Furthermore $\nu_F$ is then a symmetric measure for the perturbed process $\{Z_.^F \,
\P_x\}_{x\in E}$, which is associated with  $\mathcal E_F$ (see (\ref{3.3}.2)). For all this see
\cite{CX02} (especially Lemma 7.1 and section 2).

If in addition $F\in D(A)$, it is enough to apply Ito's formula in order to get another expression
for $Z_t^F$, namely
\begin{equation}\label{3.5}
Z_t^F \, = \,  \exp \, \{F(X_0)-F(X_t) + \, \int_0^t \, \big(AF(X_s)- \, \frac 12 \,
\Gamma(F,F)(X_s)\big) ds\} \, .
\end{equation}
 If $P_t^F$ denotes the associated ($\nu_F$ symmetric) semi-group, it holds $\nu_F$ a.s.
\begin{equation}\label{3.6}
(P_t^Fh)(x)=e^F(x) \, \E^{\P_x}\big[h(X_t) \, e^{-F(X_t)} \, M_t\big] \, ,
\end{equation}
with $$M_t=\exp \, \Big(\int_0^t \, \big(AF(X_s)- \, \frac 12 \, \Gamma(F,F)(X_s)\big) ds\Big) \,
.$$

When $\mu$ is a probability measure, $e^F\in \L^2(\nu_F)$, and a necessary condition for $\nu_F$ to
satisfy (DLSI) is thus
\begin{equation}\label{3.7}
P_t^F(e^F)=e^F \, \E^{\P_x}[M_t] \in \, \L^p(\nu_F)
\end{equation}
for all (some) $p>2$ and $t$ large enough. When $\mu$ is no more bounded one can formulate  similar
statements. For instance, if $e^F\in \L^r(\nu_F)$ for some $r>1$, then \eqref{3.7} has to hold for
some (all) $p>r$ and $t$ large enough. One can also take $r=1$ in some cases. Since the exact
formulation depends on the situation we shall not discuss it here.

A remarkable fact is that the (almost always) necessary condition \eqref{3.7} is also a sufficient
one. The next two theorems explain why. Though the proof of the first one is partly contained in
\cite{CX02} (Proposition 8.8) we shall give here the full proof for completeness.

\begin{theorem}\label{3.8}
Assume that $P_t$ is ultracontractive with $\parallel P_t\parallel_{p,\infty}=K(t,p)$ for all $p\geq
1$. Assume that H(F) is in force, $F\in D(A)$ and $M_t$ is bounded by some constant $C(t)$. Then a
sufficient condition for $\nu_F$ to satisfy (DLSI) is that $$P_t^F(e^F) \, = \, e^F \,
\E^{\P_x}[M_t] \in \, \L^q(\nu_F)$$ for some $t>0$ and some $q>2$.
\end{theorem}
\begin{proof}
Pick some $f\in \D$. Since $|f|e^{-F}\in \L^2(\mu)$ and using the Markov property, for $t>0$, $q>2$
it holds
\begin{eqnarray*}
\int \, (P_{t+s}^F(|f|))^q \, d\nu_F & = & \int \, e^{qF} \, \Big(\E^{\P_x}[M_t \, \E^{\P_{X_t}}[M_s
\, \big(e^{-F} \, |f|\big)(X'_s)]]\Big)^q \, d\nu_F \, ,\\
 & \leq & \int \, e^{qF} \, (C(s))^q \, \Big(E^{\P_x}[M_t \,
(P_s(|f| \, e^{-F}))(X_t)]\Big)^q \, d\nu_F \,  \\ & \leq & (C(s))^q \, (\parallel
P_s\parallel_{2,\infty})^q \,
\parallel f\parallel_{\L^2(\nu_F)}^q \, \int \, \big(e^F \, \E^{\P_x}[M_t]\big)^q \, d\nu_F \, .
\end{eqnarray*}
Hence
\begin{equation}\label{3.9}
\parallel P_{t+s}^F\parallel_{2,q} \, \leq \, C(s)K(s,2) \, \parallel e^F \,
\E^{\P_x}[M_t]\parallel_{\L^q(\nu_F)}\, ,
\end{equation}
and we are done.
\end{proof}

Recall that if in addition either $\mu$ is a probability measure, or $e^F\in \L^p(\nu_F)$ for some
$p>1$, condition in the Theorem is also necessary.
\medskip

When $P_t$ is only hyperbounded, the previous arguments are no more available and one has to work
harder to get the following analogue of Theorem \ref{3.8}

\begin{theorem}\label{3.10}
Assume that $P_t$ is hyperbounded. Assume that H(F) is in force, $F\in D(A)$ and that $M_t$ is
bounded by some constant $C(t)$. Assume  in addition that $e^F\in \L^r(\nu_F)$ for some $r>1$ (we
may choose $r=2$ when $\mu$ is a Probability measure).

 Then a necessary and sufficient condition for $\nu_F$ to satisfy (DLSI) is that $$P_t^F(e^F) \, = \,
 e^F \, \E^{\P_x}[M_t] \in \,
\L^p(\nu_F)$$ for some $p>2$ and some $t>0$ large enough.
\end{theorem}
\begin{proof}
The proof is based on the following elementary consequence of Girsanov theory and the variational
characterization of relative entropy (see \cite{CX02} section 2) : if $\int  f^2 \, d\nu_F=1$ and
$f$ is nonnegative, then
\begin{equation}\label{3.11}
\int \, (\sum_j \, \log \, h_j) \, f^2  d\nu_F  \leq  \frac t2 \, \mathcal E_F(f,f) + \, \log \,
\int \, f^2 \,  h_1 \,  P_t^F(h_2) \, d\nu_F \, .
\end{equation}
 Choose
$j=1,2$ , $h_1=f^{\alpha - 1}$ and $h_2=f^{\beta}$. \eqref{3.11} becomes
\begin{equation}\label{3.12}
\frac{(\alpha + \beta - 1)}{2} \, \int \,  f^2 \, \log (f^2) \, d\nu_F  \leq  \frac t2 \, \mathcal
E_F(f,f) + \, \log \, \int \, f^{1+\alpha} \, P_t^F(f^{\beta})  \, d\nu_F \, .
\end{equation}
Let $(q,s)$ a pair of conjugate real numbers. Then $$P_t^F(f^{\beta}) \, \leq \, \big(P_t^F(f^{q \,
\beta} \, e^{-\frac qs \, F})\big)^{\frac 1q} \, \big(P_t^F(e^F)\big)^{\frac 1s} \, , $$ and
accordingly
\begin{equation}\label{3.13}
\int \, f^{1+\alpha} \, P_t^F(f^{\beta})  \, d\nu_F \,  \leq \, \int \, f^{1+\alpha} \,
\big(P_t^F(f^{q \, \beta} \, e^{-\frac qs \, F})\big)^{\frac 1q} \, \big(P_t^F(e^F)\big)^{\frac 1s}
\, d\nu_F
\end{equation}
\begin{eqnarray*}
 & \leq & \Big(\int \, f^{1+\alpha} \, e^{-q \delta \, F} \, P_t^F(f^{q \, \beta} \, e^{-\frac qs
\, F}) \, d\nu_F \Big)^{\frac 1q} \, \Big(\int \, f^{1+\alpha} \,  e^{s \delta \, F} \, P_t^F(e^F)
\, d\nu_F \Big)^{\frac 1s} \\ & \leq & \Big(\int \, e^{- \, \frac{2q \delta}{1-\alpha} \, F} \,
\big(P_t^F(f^{q \, \beta} \, e^{-\frac qs \, F})\big)^{\frac{2}{1-\alpha}} \, d\nu_F \Big)^{\frac
{1-\alpha}{2q}} \, \Big(\int \,   e^{\frac{2s \delta}{1-\alpha} \, F} \,
\big(P_t^F(e^F)\big)^{\frac{2}{1-\alpha}} \, d\nu_F \Big)^{\frac {1-\alpha}{2s}}
\end{eqnarray*}
where we have used H\"{o}lder's inequality successively with $f^{1+\alpha} \, d\nu_F$ and $d\nu_F$,
and we also used $\int  f^2 \, d\nu_F=1$ to get the last expression. We have of course to choose
$\alpha < 1$. We shall also choose $\beta =1$. The first factor in the latter expression can be
rewritten
\begin{equation*}
\int \, e^{- \, \frac{2q \delta}{1-\alpha} \, F} \, \big(P_t^F(f^{q} \, e^{-\frac qs \,
F})\big)^{\frac{2}{1-\alpha}} \, d\nu_F \, = \, \int \, e^{\theta F} \, \big(\E^{\P_x}(f^q(X_t) \,
e^{-(1+\frac qs) \, F(X_t)} \, M_t)\big)^{\frac{2}{1-\alpha}} \, d\mu \, ,
\end{equation*}
with $$\theta = -  \frac{2q \delta}{1-\alpha} + \frac{2}{1-\alpha} - 2 \, .$$ Hence if we choose
$\alpha = q \delta < 1$, $\theta=0$. Furthermore $q=1+\frac qs$ and $f^q \, e^{-qF} \in \L^{\frac
2q}(\mu)$ with norm $1$, provided $q<2$. Using our hypotheses we thus obtain
\begin{equation}\label{3.14}
\int \, e^{- \, \frac{2q \delta}{1-\alpha} \, F} \, \big(P_t^F(f^{q} \, e^{-\frac qs \,
F})\big)^{\frac{2}{1-\alpha}} \, d\nu_F \, \leq \, \big( C(t) \, \parallel P_t \parallel_{\frac 2q
\, , \, \frac{2}{1-\alpha}}\big)^{\frac{2}{1-\alpha}} \, .
\end{equation}
For the second factor we choose $$\frac{2s \delta}{1-\alpha} \, < \, r \, ,$$ and since $\alpha=q
\delta$, this choice imposes $$\delta \, < \, \frac{r}{2s+rq} \quad \, \textrm{ hence } \quad \alpha
\, < \, \frac{rq}{2s+rq} \, .$$ Note that the condition $\alpha < 1$ is then automatically
satisfied. Applying H\"{o}lder again we get
\begin{equation}\label{3.15}
\int \,   e^{\frac{2s \delta}{1-\alpha} \, F} \, \big(P_t^F(e^F)\big)^{\frac{2}{1-\alpha}} \, d\nu_F
\, \leq \, \big(\int \, e^{r F} \, d\nu_F\big)^{\frac{2s\delta}{r(1-\alpha)}} \, \, \,  \big(\int \,
(P_t^F(e^F))^{p} \, d\nu_F\big)^{\frac{r(1-\alpha)-2s\delta}{r(1-\alpha)}} \, ,
\end{equation}
if $$p \, = \, \frac{2r}{r(1-\alpha)-2s\delta} \quad \textrm{ hence } \quad
\alpha=\frac{r(p-2)}{p(2(s-1)+r))} \, .$$ It remains to check that all these choices are compatible,
i.e $$\frac{r(p-2)}{p(2(s-1)+r))} \, < \, \frac{rq}{2s+rq}$$ which is easy.

Plugging \eqref{3.14} and \eqref{3.15} into \eqref{3.12} we obtain
\begin{equation}\label{3.16}
\alpha \, \int \, f^2 \, \log(f^2) \, d\nu_F \, \leq \, t \, \mathcal E_F(f,f) \, + \, 2A \, ,
\end{equation}
where $$A = \frac 1q \, \log\big(C(t) \, \parallel P_t \parallel_{\frac 2q \, , \,
\frac{2}{1-\alpha}}\big)+ \, \frac{\alpha}{q} \, \log\big(\parallel e^F\parallel_{\L^r(\nu_F)}\big)
+ \frac 1s \,  \log\big(\parallel P_t^F(e^F)\parallel_{\L^p(\nu_F)}\big) \, .$$ For a fixed $p$ we
may choose any pair $(q,s)$ with $q<2$, and the corresponding $\alpha$ yields the result for $$t \,
\geq \, \frac a2 \, \log\Big(\frac{q(1+\alpha)}{(2-q)(1-\alpha)}\Big) \, ,$$ according to Gross
theorem, if $\mu$ satisfies (DLSI) with constants $(a,b)$.
\end{proof}
\medskip

\begin{remark}\label{3.17}
Unfortunately the previous methods cannot furnish the best constants. In particular we cannot get
(TLSI) even when $\mu$ satisfies (TLSI).
\end{remark}

In view of the previous remark it is thus natural to look at the spectral gap properties too. The
final result we shall recall is Lemma 2.2 in \cite{Aid01}.

\begin{theorem}\label{3.18}
Assume that $\mu$ is a probability measure satisfying (SGP).Assume that H(F) is in force and
$\Gamma(F,F) \in \L^1(\mu)$. Then $\nu_F$ satisfies (WSGP).
\end{theorem}

One can use Theorems \ref{3.8} (or \ref{3.10}) and \ref{3.18} together in order to show that the
general Boltzmann measure satisfies (TLSI) provided $\mu$ is a Probability measure. Otherwise one
has to consider various reference measures $\mu$, as it will be clear in the next sections.
\bigskip

\section{\bf The ``Well Method''.}\label{IV}

Our aim in this section is to get sufficient general conditions for \eqref{3.7} to hold. To this end
we shall slightly modify the ``Well Method'' of \cite{KKR}, i.e. use the martingale property of the
Girsanov density. In the sequel we assume that $F\in D(A)$ satisfies H(F).
\smallskip

The main assumption we shall make is the following, for all $x$
\begin{equation}\label{4.1}
\frac 12 \, \Gamma(F,F)(x) \, - \, AF(x) \, \geq \, -c \, > - \infty \, .
\end{equation}
It follows that $M_t \, \leq \, e^{ct} \, = \, C(t)$.

Now we define $\lambda(x)$ by the relation , $$\frac 12 \, \Gamma(F,F)(x) \, - \, AF(x) \, = \,
\lambda(x) \, F(x)  \, .$$ Note that if $F(x)\leq 0$, $P_t^F(e^F)(x) \, \leq \, C(t)$ so that the
contribution of the $x's$ with $F(x)\leq 0$ belongs to $\L^{\infty}(\nu_F)$. So we may and will
assume that $F(x)>0$.
\medskip

For $0<\varepsilon<1 $ define the stopping time $\tau_x$ as

\begin{equation}\label{4.2}
\tau_x \, = \, \inf \{ \, s>0 \, , \, (\frac 12 \, \Gamma(F,F) \, - \, AF)(X_s) \, \leq \,
\varepsilon \, \lambda(x) \, F(x) \, \textrm{ or } \, F(X_s) \, \leq \, \varepsilon \, F(x)\} \, .
\end{equation}

First we assume that $(\frac 12 \, \Gamma(F,F) \, - \, AF)(x) \, > \, 0$ . In this case $\tau_x > 0
\, $ $\, \P_x$ a.s. Introducing the previous stopping time we get
\begin{equation*}
\E^{\P_x}[M_t]  =  \E^{\P_x}[M_t \, \BBone_{t<\tau_x}]+\E^{\P_x}[M_t \, \BBone_{\tau_x \leq t}] = A
+ B \, ,
\end{equation*}
with
\begin{equation}\label{4.4}
A \, = \, \E^{\P_x}[M_t \, \BBone_{t<\tau_x}] \, \leq \, \exp \, -\big( \varepsilon \, t \,
\lambda(x) \, F(x)\big)  \, ,
\end{equation}
and
\begin{equation}\label{4.5}
B \, =  \, \E^{\P_x}[M_t \, \BBone_{\tau_x \leq t}]
\end{equation}
\begin{eqnarray*}
 \leq & e^{ct} \,  \E^{\P_x}[\exp \,
\Big(\int_0^t \, \big(AF \, - \, \frac 12 \Gamma(F,F) \, + c \big)(X_s) \, ds\Big) \, \BBone_{\tau_x
\leq t}] \nonumber \\ \leq & e^{ct} \, \E^{\P_x}[\exp \, \Big(\int_0^{\tau_x} \, \big(AF \, - \,
\frac 12 \Gamma(F,F) \, + c \big)(X_s) \, ds\Big) \, \BBone_{\tau_x \leq t}] \nonumber
\\ \leq & e^{ct} \, \E^{\P_x}[\exp \, \Big(\int_0^{\tau_x} \, \big(AF \, - \, \frac 12
\Gamma(F,F)\big)(X_s) \, ds\Big) \, \BBone_{\tau_x \leq t}] \\  = & e^{ct} \, \E^{\P_x}[M_{\tau_x}
\,  \BBone_{\tau_x \leq t}].
\end{eqnarray*}

But $e^{-F(X_s)} \, M_s$ is a  $\L^2$ (thanks to \ref{4.1}) $\P_x$ martingale. Hence, according to
Doob's Optional Stopping Theorem
\begin{equation}\label{4.6}
\E^{\P_x}[e^{-F(X_{\tau_x})} \, M_{\tau_x} \, \BBone_{\tau_x \leq t}] \, \leq \,
\E^{\P_x}[e^{-F(X_{t\wedge \tau_x})} \, M_{t\wedge \tau_x}] \, = \, e^{-F(x)} \, .
\end{equation}

According to \eqref{4.2}, $$e^{-F(X_{\tau_x})} \, \geq \, e^{- \, \varepsilon \, F(x)} \, ,
$$ so that thanks to \eqref{4.6}, $$\E^{\P_x}[M_{\tau_x} \,  \BBone_{\tau_x \leq t}] \, \leq \,
e^{-(1-\varepsilon) F(x)}  \, .$$ Using this estimate in \eqref{4.5} and using \eqref{4.4} we
finally obtain
\begin{equation}\label{4.7}
\E^{\P_x}[M_t] \, \leq \, e^{- \, \varepsilon \, t \, \lambda(x) \, F(x)} \, + \,  e^{ct} \,
e^{-(1-\varepsilon) F(x)} \, .
\end{equation}

Finally if $(\frac 12 \, \Gamma(F,F) \, - \, AF)(x) \, < \, 0$ we certainly have
\begin{equation*}
\E^{\P_x}[M_t] \, \leq \, e^{ct} \, e^{- \, \varepsilon \, t \, \lambda(x) \, F(x)} \, ,
\end{equation*}
since in this case $\lambda(x)<0$ while we assume $F(x)>0$.
\smallskip

We have thus obtained choosing first $\varepsilon=r/p$,
\smallskip

\begin{theorem}\label{4.8}
Assume that H(F) and \eqref{4.1} are fulfilled. Assume in addition that there exists some $0<r$ such
that $e^{F} \in \L^r(\nu_F)$. Then $e^{F} \, \E^{\P_x}[M_t]  \in  \L^p(\nu_F)$ as soon as $$\int \,
e^{(p-2)F} \, e^{-(rt/p) \, (\frac 12 \, \Gamma(F,F) \, - \, AF)} \, d\mu \, < \, +\infty \, .$$ In
particular $\nu_F$ satisfies (DLSI) as soon as $$\int \, e^{\beta F} \, e^{- \, \lambda \, (\frac 12
\, \Gamma(F,F) \, - \, AF)} \, d\mu \, < \, +\infty \, ,$$ for some $\beta>0$ and some $\lambda>0$.
Furthermore if the previous holds for all pair $(\beta,\lambda)$ of positive real numbers, then
$P_t^F$ is immediately hyperbounded (i.e. $P_t^F$ is bounded from $\L^2(\nu_F)$ in $\L^p(\nu_F)$ for
all $t>0$ and all $p>2$).
\end{theorem}
\medskip

\begin{remark}\label{4.9}
This result extends previous ones obtained by Davies \cite{Dav} (especially Theorem 4.7.1 therein)
in the ultracontractive context, by Rosen \cite{Ro76} in the hyperbounded context (based on deep
Sobolev inequalities available in $\R^N$) or by Kusuoka and Stroock \cite{KS85}. In addition it is
an ``almost'' necessary condition too, in the sense of the next result.
\end{remark}
\bigskip

\begin{theorem}\label{4.10}
Assume that H(F) holds and that there exists some $1<r$ such that $e^{F} \in \L^r(\nu_F)$. A
necessary condition for $\nu_F$ to satisfy (DLSI) is $$\int \, e^{\beta F(x)} \, e^{- \, \lambda \,
(\frac 12 \, \Gamma(F,F) \, - \, AF)(x)} \, \P_x^{2+\beta}(\tau_x>\lambda/2(2+\beta)) \, d\mu \, <
\, +\infty \, ,$$ for some $\beta>0$ and some $\lambda>0$, where $\tau_x$ is the stopping time
defined by
$$\tau_x \, = \, \inf \, \{s\geq 0 \, s.t. \, (\frac 12 \, \Gamma(F,F) \, - \, AF)(X_s) \, \geq \, 2
\lambda(x) \, F(x)\} \, , $$ $\lambda(x)$ being defined as $$\frac 12 \, \Gamma(F,F)(x) \, - \,
AF(x) \, = \, \lambda(x) \, F(x) \, .$$
\end{theorem}

\begin{proof}
It is enough to remark that $\tau_x=0$ if $\lambda(x)\leq 0$ and then write for $\lambda(x)>0$
\begin{eqnarray*}
\E_x[M_t] & \geq & \E_x[M_t \, \BBone_{t<\tau_x}]\\ & \geq & e^{-2t \lambda(x) F(x)} \,
\E_x[\BBone_{t<\tau_x}] \, ,
\end{eqnarray*}
and then to apply the necessary part of Theorem \ref{3.10}.

\end{proof}
\bigskip

\section{\bf A direct approach for the sufficient condition and others consequences.}\label{ls}

In the previous two sections we used the hyperbounded point of view. As suggested by an anonymous
referee, Theorem \ref{4.8} can be directly obtained by using logarithmic Sobolev inequalities.

Indeed assume that
\begin{equation}\label{ls1}
\int \, f^2 \, \log f^2 \, d\mu \, \leq \, C_1 \, \int \, \Gamma(f,f) \, d\mu \, + \, C_2 \, ,
\end{equation}
for all nice $f$ such that $\int \, f^2 \, d\mu \, = \, 1 $ . Take $f \, = \, e^{-F} \, g$ for some
$g$ such that $\int \, g^2 \, d\nu_F \, = \, 1$ . Thanks to the chain rule, i.e. $$\int \,
\varphi'(f) \, Af \, + \, \frac 12 \, \varphi''(f) \, \Gamma(f,f) \, d\mu \, = \, 0 \, $$ it is easy
to see that \eqref{ls1} can be rewritten

\begin{equation}\label{ls2}
\int \, g^2 \, \log g^2 \, d\nu_F \, \leq \, C_1 \, \int \, \Gamma(g,g) \, d\nu_F \, + \, \int \,
g^2 \, \left(2C_1\big(AF \, - \, \frac 12 \, \Gamma(F,F)\big) \, + \, 2F\right) \, d\nu_F \, + \,
C_2 \, .
\end{equation}

Introducing some $0 < \varepsilon < 1$, we write the second integral in the right hand side
$$\varepsilon \, \int \, g^2 \, \frac{1}{\varepsilon} H \, d\nu_F \, ,$$ and use Young's inequality
in order to get

\begin{equation}\label{ls3}
(1-\varepsilon) \, \int \, g^2 \, \log g^2 \, d\nu_F \, \leq \,
\end{equation}
\begin{eqnarray*}
& \leq & C_1 \, \int \, \Gamma(g,g) \, d\nu_F \, + \, \varepsilon \, e^{-1} \, \int \, e^{\frac{2
C_1}{\varepsilon}\big(AF \, - \, \frac 12 \, \Gamma(F,F)\big)+2(\frac 1{\varepsilon} - 1)F} \, d\mu
\, + \, C_2 \, ,
\end{eqnarray*}
and we recover Theorem \ref{4.8} since we may choose $\varepsilon$ arbitrarily close to $1$ and
independently $C_1$ arbitrarily large. Actually in Theorem \ref{4.8}, since \eqref{4.1} is
fulfilled, we may choose any $\lambda' > \lambda$. The only difference here is that we do not need
to assume \eqref{4.1}, but in contrast, we have to assume that $\lambda$ is large enough.

The above proof is given with less details than the previous martingale proof. Actually both are
short and elementary. The main advantage of the martingale point of view is to indicate how to get a
necessary condition.
\medskip

However it is interesting at this point to compare our condition for (DLSI) and known results on
(SGP) obtained by Gong and Wu \cite{GW00} for Feynman-Kac semigroups.

The unitary transform $U \, : \, \L^2(E,d\mu) \, \rightarrow \, \L^2(E,d\nu_F)$ defined by $U(f)=e^F
\, f$ satisfies $$\int \, \Gamma(U(f),U(g)) \, d\nu_F \, = \, \int \, \big(\Gamma(f,g) \, + \, V_F
\, fg\big) \, d\mu$$ where $V_F= \Gamma(F,F) \, - \, 2 \, AF$ . The latter Dirichlet form is the one
associated with the Schr\"{o}dinger operator $H_F=\, A \, + \, V_F$. Since $U$ is unitary the
spectrum of $H_F$ on $\L^2(d\mu)$ and the one of $-A_F$ on $\L^2(\nu_F)$ coincide. Hence the
existence of a spectral gap for $\nu_F$ follows from Corollary 6 in \cite{GW00}, namely

\begin{proposition}\label{ls4}
Let $\mu$ be a probability measure satisfying (TLSI) (i.e. \eqref{ls1} with $C_2=0$) and assume that
H(F) holds. If
$$ \int \, e^{(2 \, C_1+\varepsilon)(\frac 12 \Gamma(F,F) \, - \, AF)^{-}} \, d\mu \, < \, +\infty
$$ for some $\varepsilon > 0$ then $\nu_F$ satisfies (SGP). This result holds in particular when
\eqref{4.1} is satisfied.
\end{proposition}
It follows in particular that, provided $F$ is bounded below, the condition in Proposition \ref{ls4}
is implied by the condition in Theorem \ref{4.8} without assuming \eqref{4.1}, but assuming that
$\lambda > 2 C_1$.

\begin{corollary}\label{ls5}
If $\mu$ satisfies (TLSI) (or equivalently $P_t$ is hypercontractive) and H(F) holds, then $$\int \,
e^{\beta F} \, e^{- \, \lambda \, (\frac 12 \, \Gamma(F,F) \, - \, AF)} \, d\mu \, < \, +\infty \,
,$$ for some $\beta>0$ and $\lambda>0$ is a sufficient condition for $\nu_F$ to satisfy (TLSI)
provided in addition
\begin{enumerate}
\item[(1)] \quad either $\frac 12 \, \Gamma(F,F) \, - \, AF$ is bounded from below and $e^F\in
\L^r(\nu_F)$ for some $r>0$ ,

\item[(2)] \quad or $F$ is bounded below and $\lambda > 2 C_1$ where $C_1$ is the optimal constant
in (TLSI) for $\mu$.
\end{enumerate}
\end{corollary}
 The interested reader will find a stronger statement (Theorem 5) in \cite{GW00}, but
with less tractable hypotheses.
\medskip

\bigskip

\section{\bf Examples: $\R^N$ valued Boltzmann measures. }\label{V}

In this section we shall deal with the $\R^N$ valued case, i.e. $E=\R^N$, $dx$ is Lebesgue measure,
$A=\frac 12 \, \Delta$ is one half of the Laplace operator and $\nabla$ is the usual gradient
operator. $\P_x$ is thus the law of the Brownian motion starting at $x$, whose associated semigroup
$P_t$ is $dx$ symmetric and ultracontractive with $\parallel P_t\parallel_{2,+\infty}=(4\pi \,
t)^{-\frac N4}$ . $\D$ is the algebra generated by the usual set of test functions and the
constants.
\smallskip

(TLSI) can thus be written
\begin{equation*}
\int \, f^2 \, \log \Big(\frac{f^2}{\parallel f\parallel_{\L^2(\nu_F)}^2}\Big) \, e^{-2F} \, dx \,
\leq \, a \, \int \, |\nabla f|^2 \, e^{-2F} \, dx \, .
\end{equation*}

Note that Lebesgue measure satisfies a family of logarithmic Sobolev inequalities i.e. for all $\eta
> 0$ and all $f$ belonging to $\L^1(dx)\cap \L^{\infty}(dx)$ such that $\int \, f^2 \, dx \, = \, 1$
\begin{equation*}
\int \, f^2 \, \log f^2 \, dx \, \leq \, 2\eta \, \int \, |\nabla f|^2 dx \, + \, \frac N2 \,
\log\left(\frac1{4\pi \eta}\right) \, ,
\end{equation*}
see e.g. \cite{Dav} Theorem 2.2.3.

 In the sequel we will consider functions $F$ that are of class $C^2$ and according to
Proposition \ref{3.2} we shall then (if necessary) add to $F$ some bounded perturbation. Furthermore
in this particular finite dimensional situation we may replace H(F) by the following Lyapounov
control:
\begin{equation}\label{5.1}
\textrm{there exists some } \psi \textrm{ such that } \psi(x) \, \rightarrow \, +\infty \, \textrm{
as } \, |x| \, \rightarrow \, +\infty \, \, ,
\end{equation}
$$\textrm{ and } \, \Delta \psi(x) \, - \, (\nabla F \, . \, \nabla \psi)(x) \, \leq \, K \, < \,
+\infty\, \textrm{ for all } \, x \, .$$
\smallskip

In order to complete the picture, we have to describe some sufficient conditions allowing to tight
the logarithmic Sobolev inequality.

One is given by Theorem \ref{3.18}. Indeed if $\nu_U(dx)=e^{-2U(x)}dx$ is another Boltzmann measure
satisfying (SGP) and $$\int \, |\nabla U|^2 \, d\nu_U \, < \, +\infty$$ a sufficient condition for
$\nu_F$ to satisfy (WSGP) is
\begin{equation*}
\int \, |\nabla F|^2 \, d\nu_U \, < \, +\infty \, ,
\end{equation*}
since $d\nu_F=e^{-2(F-U)}d\nu_U$. It is thus not difficult to guess that (WSGP) holds for \emph{any}
Boltzmann measure (such that $F$ is smooth). This result is actually true and shown (using another
route) in \cite{RW01} Theorem 3.1 and Remark (1) following this theorem. Hence

\begin{proposition}\label{5.2}
For a Boltzmann measure $\nu_F$ with $F\in C^2$ , (WSGP) is satisfied. Consequently (DLSI) and
(TLSI) are equivalent.
\end{proposition}

It is nevertheless interesting, at least for counter examples to know some sufficient conditions for
the usual (SGP). If $N=1$ a necessary and sufficient condition was obtained by Muckenhoupt (see
\cite{Ane} chapter 6). We recall below a tractable version due to Malrieu and Roberto of this result
as well as its $N$ dimensional counterpart

\begin{proposition}\label{5.3}
Let $F$ of $C^2$ class.
\begin{enumerate}
\item[(1)] \quad (see \cite{Ane} Theorem 6.4.3) \quad If $N=1$, $|F'(x)|>0$ for $|x|$ large enough
and $\frac{F''(x)}{|F'(x)|^2}$ goes to $0$ as $|x|$ goes to $\infty$, then $\nu_F$ satisfies (SGP)
if and only if $$\liminf_{|x| \, \rightarrow \, +\infty} |F'(x)|^2 \, = \, C \, > \, 0 \, .$$
\item[(2)] \quad
 (see e.g. \cite{Ku02} Proposition 3.7) \quad For any $N$,
 if $$\liminf_{|x| \, \rightarrow \, +\infty} (|\nabla \, F|^2 \, - \,
\Delta \, F) \, = \, C \, > \, 0 \, ,$$ then (SGP) holds for $\nu_F$.
\end{enumerate}
\end{proposition}

Now if we want to use Proposition \ref{ls4} we may choose $d\mu=(1/Z_{\rho}) \, e^{-2\rho |x|^2} \,
dx$ which is known to satisfy (TLSI) with constant $C_1=1/2\rho$, and is associated to the generator
$$A_{\rho} \, = \, \frac 12 \, \Delta \, - \, 2\rho \, x.\nabla \, .$$ We thus have to look at
\begin{equation}\label{sg}
\frac 12 \, |\nabla (F-\rho |x|^2)|^2 \, - \, A_{\rho}(F-\rho |x|^2) \, = \, \frac 12 \,
\big(|\nabla F|^2 \, - \, \Delta F\big) \, - 2\rho^2 |x|^2 \, + \, \rho N \, .
\end{equation}
Note thus that we cannot recover \ref{5.3}(2).
\medskip

According to Proposition \ref{5.2} and the previous sections we know that a sufficient condition for
(TLSI) to hold is the integral condition in Theorem \ref{4.8} (assuming in addition one of
conditions (1) and (2) in Corollary \ref{ls5}), while a necessary one is given in Theorem \ref{4.9}.
Up to our knowledge, except the bounded perturbation recalled in Proposition \ref{3.2}, three others
family of sufficient conditions have been given for $\nu_F$:
\begin{itemize}
\item the renowned Bakry-Emery criterion saying that (TLSI) holds as soon as $F$ is uniformly
convex, i.e. $Hess(F) \, \geq \, K \, Id$ for some $K>0$, \item Wang's results (see \cite{Wa01}
Theorem 1.1 for this final version) saying that provided $Hess(F) \, \geq \, -K \, Id$ for some
$K\geq 0$, a sufficient condition is $$\int \, e^{\varepsilon |x|^2} \, d\nu_F \, \leq \, +\infty$$
for some $\varepsilon \, > \, K$ , \item the beautiful Bobkov-G\"{o}tze criterion for $N=1$, and its
weak version due to Malrieu and Roberto (see \cite{Ane} Theorem 6.4.3) saying that if $|F'(x)|>0$
for $|x|$ large enough and $\frac{F''(x)}{|F'(x)|^2}$ goes to $0$ as $|x|$ goes to $\infty$, then
$\nu_F$ satisfies (TLSI) if and only if there exists some $A$ such that $$\frac{F}{|F'|^2} \, + \,
\frac{\log|F'|}{|F'|^2} $$ is bounded on $\{|x| \, \geq \, A\}$.
\end{itemize}
\medskip

It is not difficult to see that our results contain Malrieu-Roberto result.

It is also easy to see that if $Hess(F)(x) \, \geq \, \rho \, Id$ for some positive $\rho$ and all
$x$, then $$|\nabla F|^2(x) \, \geq \, 2\rho F(x) \, - \, C \, ,$$ for some constant $C$. Hence if
$F$ is uniformly convex and such that $$|\Delta F|(x) \, \leq \, (1-\varepsilon) \, |\nabla F|^2(x)
\, + \, c(F) \, ,$$ for some $\varepsilon > 0$ , all $x$ , and some constant $c(F)$, we recover the
result by Bakry-Emery. The same holds for Wang's result if the perturbed $F+\frac 12 \,
(K+\varepsilon) |x|^2$ is a nice uniformly convex function as before.

Unfortunately, it is not difficult to build uniformly convex functions such that
$$\limsup_{|x|\rightarrow +\infty} \, \left(\frac {\Delta F}{|\nabla F|^2}\right) \, = \, +\infty \,
.$$ Actually the counter examples built by Wang are such that the previous property holds.
\medskip

\begin{remark}\label{5.4}

Assume that $\lim_{|x|\to +\infty} \, F(x) \, = \, +\infty$ .

Applying Theorem \ref{4.8} we see that $P_t^F$ is hypercontractive in particular as soon as
$$|\nabla F|^2(x) \, - \, \Delta F(x) \, \geq \, \eta \, F(x) \, - \, c \, ,$$ for some constant $c$
and some $\eta
>0$.

As we remarked in Theorem \ref{4.8}, we can get conditions for \emph{immediate} hypercontractivity,
for instance $P_t^F$ will be \emph{immediately} hypercontractive as soon as
$$|\nabla F|^2(x) \, - \, \Delta F(x) \, \geq \, G(F(x))  \, ,$$ for some  function $G$ such that
$$\lim_{y\to +\infty} \, \frac{G(y)}{y} \, = \, +\infty \, .$$ One can also see from \cite{KKR} that
a condition like $$\int^{+\infty} \, \frac{y}{G(y) \, g'(g^{-1}(y))} \, dy \, < \, +\infty \quad
\textrm{ for some $g$ satisfying } \quad \int^{+\infty} \, e^{- \, g(y)} \, dy \, < \, +\infty \,
,$$ for the function $G$ we have introduced above, implies that $P_t^F$ is \emph{ultracontractive}.
This result with $G(y)=y^{\theta}$ for some $\theta>1$ (take then $g(y)=e^y$) is contained in
\cite{Dav} Theorem 4.7.1.

As shown in \cite{BCR} the same control but with $0 <\theta < 1$ yields a weaker form of
hypercontractivity.
\end{remark}
\bigskip

As we discussed before, if our results can only be partly compared (at least easily) with existing
ones in the bounded below curvature case (i.e. when the Hessian is bounded from below), they allow
to look at interesting examples in the unbounded curvature case. We shall below discuss such a
family of examples. But first we recall a basic estimate for the Brownian motion that allows us to
give a precise meaning to the necessary condition stated in Theorem \ref{4.10}.
\smallskip

\begin{lemma}\label{5.5}

For a standard Brownian motion $B_s$ on $\R^N$, there exists a constant $\theta_N$ such that $$\P \,
(\sup_{0\leq s \leq t} \, |B_s| \, < \, A) \, \geq \, e^{-\theta_N \, \frac{t}{A^2}} \, .$$
\end{lemma}

\medskip

\begin{example}\label{5.6}
Let us consider on $\R^+$ the potential $F_{\beta}(x) \, = \, x^2 + \, \beta \, x \, \sin(x)$
extended by symmetry to the full real line. We shall only look at its behaviour near $+\infty$.

The derivatives are given by $F'_{\beta}(x) \, = \, (2+\beta \, \cos(x)) \, x + \, \beta \, \sin(x)$
and $F''_{\beta}(x) \, = \, -\beta \, x \, \sin(x) \, + \, 2(1+\beta \, \cos(x))$ . Hence $-\infty
\, = \, \liminf_{x\rightarrow +\infty} F''(x)$.

For $|\beta| < 2$ we may apply Malrieu-Roberto result (or Theorem \ref{4.8}) and show that (TLSI)
holds.

For $|\beta| \, \geq \, 2$ the hypotheses of Theorem \ref{4.8} are no more satisfied. Indeed
$${F'}^2(x) - F''(x) \, = \, (2+\beta \cos(x))^2 \, x^2 \, + \, (4+2\beta \cos(x) - \beta \sin(x)) x +
\,  h(x)$$ where $h$ is bounded, can be very negative for the $x$'s such that $2+\beta \cos(x) =0$.
\smallskip

We shall discuss below the case $\beta=-2$ in details. Instead of using Theorem \ref{4.10} we shall
directly study $P_t^F(e^F)$ for $F=F_{-2}$.

 Introduce $x_k=2k\pi$. Then for $k$ large enough one can
find $\varepsilon$ small enough and some constant $c$ such that
\begin{eqnarray}\label{5.160}
\textrm{ for all $y$ such that }  1/2 \, k^{-\frac 12} \leq \, y-x_k \, \leq 3/2 \, k^{-\frac 12} \,
\textrm{ it holds }\\ F''(y)\geq (1-\varepsilon) k^{\frac 12} \quad \textrm{ and } \quad |F'(y)|\leq
c \, .\nonumber
\end{eqnarray}
Introduce the stopping times $\tau_k = \inf \, \{s\geq 0 \, , \, |X_s - y|\geq \frac 14 \, k^{-\frac
12}\}$ . Then according to \eqref{5.160}, for $3/4 \, k^{-\frac 12} \leq \, y-x_k \, \leq 5/4 \,
k^{-\frac 12}$
\begin{eqnarray*}
\E^{\P_{y}}(M_t) & \geq & \E^{\P_{y}}(M_t \, \BBone_{t < \tau_k}) \\ & \geq & e^{\frac t2 \,
((1-\varepsilon) k^{\frac12} - c^2)} \, \P_{y}(t < \tau_k) \\ & \geq & e^{\frac t2 \,
((1-\varepsilon) k^{\frac12} - c^2)} \, e^{-\, 4 \theta t k}
\end{eqnarray*}
for the constant $\theta$ appearing in Lemma \ref{5.5}. It follows
\begin{equation*}
\int^{+\infty} \, e^{(q-2)F}\big(\E^{\P_x}[M_t]\big)^q \, dx \, \geq  \, \frac 12 \, \sum_k \,
k^{-\frac 12} \, e^{4\pi^2 \, (q-2-\varepsilon) \, k^2} \, e^{- 4 q \theta t k} \, = \, +\infty \, .
\end{equation*}
Hence (DLSI) does not hold.
\smallskip

For $\beta=2$ the discussion is similar, while for $|\beta|>2$ it is a little bit different. Indeed
(again with $\beta<0$) this time if $F'(x_k)=0$, on  $2 k^{-\frac 34}\leq y-x_k \leq k^{-\frac 34}$
we have $F''(y)\geq (1-\varepsilon) k$ while $|F'(x)|\leq c \, k^{\frac 14}$. Hence we can prove as
before that $\E^{\P_{y}}(M_t) \, \geq \, C \, e^{-c' \theta t \, k^{\frac 32}}$ for some constants
$C$ and $c'$ and conclude again that (DLSI) does not hold.
\medskip

\end{example}

\bigskip
\bibliographystyle{plain}

\end{document}